 \documentclass[a4paper, 10pt ]{article}
\usepackage{graphicx}
\usepackage{amsmath} 
\usepackage{subfig}

\def\bfu{\mbox{\boldmath$u$}}
\def\bfg{\mbox{\boldmath$g$}}
\def\J{\mbox{\boldmath $J$}}
\def\I{\mbox{\boldmath $I$}}
\def\n{\mbox{\boldmath $n$}}
\def\R{\mbox{\boldmath $R$}}

\newtheorem{theorem}{Theorem}[section]

\newtheorem{definition}[theorem]{Definition}

\title{Toward a mathematical analysis for a model of    
suspension flowing down an inclined plane}

\author{Kaname Matsue\thanks{Institute of Mathematics for Industry / International Institute for Carbon-Neutral Energy Research (WPI-I$^2$CNER), Kyushu University, Fukuoka 819-0395, Japan, (kmatsue@imi.kyushu-u.ac.jp)}
\and Kyoko Tomoeda\thanks{Institute for Fundamental Sciences, Setsunan University, Osaka
572-8508, Japan, (tomoeda@mpg.setsunan.ac.jp)}}

\begin{document}


\maketitle

\begin{abstract}
We consider the Riemann problem of the dilute approximation equations with spatiotemporally  
dependent volume fractions from the full model of suspension, in which the particles settle to 
the solid substrate and the clear liquid film flows over the sediment 
[Murisic et al., J. Fluid. Mech. {\bf{717}}, 203--231 (2013)]. We present a 
method to find shock waves, rarefaction waves for the Riemann problem of this 
system. Our method is mainly based on [Smoller, Springer-Verlag, New York, second 
edition, (1994)]. 
\end{abstract}


\pagestyle{myheadings}
\thispagestyle{plain}

\section{Introduction}
\label{1}
We are concerned here with the two dimensional motion of a suspension flowing down an inclined plane under the effect of gravity. 
To describe the problem we choose a coordinate system $(x, y)$, where the $x$-axis is along a plane with a inclination angle $\alpha$ 
$\left(0<\alpha<\frac{\pi}{2}\right)$ and the $y$-axis is perpendicular to this plane.
The motion of suspension is governed by the following partial differential equations 
\begin{eqnarray}
&& \nabla p - \nabla\cdot[\mu(\phi)(\nabla \bfu + \nabla \bfu^{\top})]=
\rho(\phi)\bfg, \nonumber\\
&& \partial_t\phi + \bfu\cdot\nabla\phi + \nabla\cdot\J =0, \label{NS}\\
&& \nabla\cdot\bfu=0, \hspace{0.5cm}\hbox{ in }\hspace{0.3cm}
0 < y < h(x, t), \hspace{0.3cm}t\geq 0.  \nonumber 
\end{eqnarray}
Here $\bfu=(u, v)^{\top}$ is the volume averaged velocity and $p$ is the pressure of fluid and $h(x, t)$ is the total suspension thickness. 
$\phi$ is the particle volume fraction and $\J=(\J_x, \J_y)^{\top}$ is the particle flux and  
$\bfg=g(\sin\alpha, -\cos\alpha)^{\top}$ is the acceleration of gravity. 
$\mu(\phi)$ is the viscosity of fluid and $\rho(\phi)=\rho_{p}\phi+ \rho_{f}(1-\phi)$, where $\rho_{f}$ and $\rho_{p}$ are the density of fluid and particles respectively. 
The boundary condition on the wall is the non-slip and no-penetration condition   
\begin{eqnarray}
\bfu=(0, 0)^{\top},\hspace{0.5cm}\hbox{ at }\hspace{0.3cm}y=0.  \label{nonslip}
\end{eqnarray}
The dynamical and kinematic conditions on the free surface are
\begin{eqnarray}
\left(-p\I + \mu(\phi)(\nabla\bfu+ \nabla\bfu^{\top})\right)\n &=& 0,
\hspace{0.5cm}\hbox{ at }\hspace{0.3cm}y=h(x,t),  \label{balance}\\
\partial_t h + u\partial_x h -v &=& 0, 
\hspace{0.5cm}\hbox{ at }\hspace{0.3cm}y=h(x,t),  \nonumber \label{kine}
\end{eqnarray}
where $\I$ is the identity matrix and $\n$ is the outward unit normal vectors 
to the free surface. For the particle fluxes, the no-flux boundary conditions 
at the wall and free surface are also imposed : 
\begin{eqnarray}
J\cdot\n = 0, \hspace{0.5cm}\hbox{ at }\hspace{0.3cm}y=0
\hspace{0.1cm}\hbox{ and }\hspace{0.1cm}y=h(x,t).  \label{noflux}
\end{eqnarray}
\par 
To explain the mechanisms of suspensions, some approximation 
equations are derived from the full model {\eqref{NS}}--{\eqref{noflux}}. 
Murisic et al. {\cite{MPPB}} derived the dilute approximation equation which is the system of conservation laws : 
\begin{eqnarray}
\partial_t h + \partial_x \bigg(\displaystyle\frac{1}{3}h^3\bigg) &=& 0, 
\label{CL} \\
\partial_t n + \partial_x \bigg(\displaystyle\sqrt{\frac{2}{9C}}(nh)^{3/2}\bigg)&=& 0, \label{PV}   
\end{eqnarray}
where $C=\frac{2(\rho_p-\rho_f)\cot\alpha}{9(\rho_pK_c)}$ is the buoyancy parameter and 
$K_c$ is constant and $n=\phi h$. This dilute approximation equation focuses on the settled  
regime in which particles settle to the solid substrate and
the clear liquid film flows over the sediment. 
In {\cite{MPPB}}, the authors solved {\eqref{CL}} exactly with the 
initial data $h(x, 0)=1$ for $0\leq x\leq 1$, $h(x, 0)=0$ otherwise, 
and the exact solution for $h$ is given by 
\begin{eqnarray*}
h(t,x) = 
\begin{cases}
1         &\hspace{0.5cm}t\leq x \leq x_{\ell}, \\
\displaystyle\sqrt{\frac{x}{t}}
          &\hspace{0.5cm}0< x < \min (t, x_{\ell}),  \label{solution}\\
0         &\hspace{0.5cm}\hbox{otherwise}, 
\end{cases}
\end{eqnarray*} 
for $t\geq 0$, where $x_{\ell}$ denote the liquid front position which is given by 
$x_{\ell}=1+\frac{t}{3}$ for $0\leq t\leq \frac{3}{2}$, 
$x_{\ell}=\left(\frac{9t}{4}\right)^{1/3}$ for $\frac{3}{2}<t$.    
One of the earlier examples for solution {\eqref{solution}} is given by Huppert {\cite{H1}} 
for the flow of a constant volume of viscous fluid down a constant slope.  
The authors in {\cite{MPPB}} also obtain the exact solution $n$ of {\eqref{PV}} with the initial data 
$n(x, 0) = f_0 h(x, 0)$ and some given value $f_0\ll 1$. \par 
Our aim in this paper is to cover the solution of the system {\eqref{CL}}--{\eqref{PV}} 
when the initial volume fraction $\phi(0,x)$ is a variable satisfying $0<\phi<1$.  
For this system, only exact solutions obtained for the fixed initial volume fraction $\phi(x, 0)=f_0$ are treated in {\cite{MPPB}}.    
On the other hand, in mathematical theory, it is known that the general $m\times m$ system of the hyperbolic conservation laws   
\begin{eqnarray*}
\partial_tU + \partial_x(F(U)) = 0  
\end{eqnarray*}
has a discontinuous solution such as a shock wave and a smooth solution such as a rarefaction wave, 
where $U= (U_1,\cdots, U_m)^{\top}\in\R^m$, $(x,t)\in \R\times \R_{+}$   
and $F(U) = (F_1(U),\cdots, F_m(U))^{\top}$ is a vector-valued function which is $C^2$ in some open subset $D\subset \R^m$ 
(see {\cite{La1}}, {\cite{S}}). In order to cover the solution of the system {\eqref{CL}}--{\eqref{PV}}, 
we consider the case where the solutions have a discontinuity, and hence we deal with the weak solution of the system 
which is defined by {\eqref{weak}} below. 
Applying mathematical theories established in {\cite{La1}}, {\cite{S}} to the system {\eqref{CL}}--{\eqref{PV}}, 
we give a construction method of weak solutions consisting of simple waves such as shock waves and rarefaction waves.
\par 
The organization of this paper is as follows. In Section {\ref{2}}, we formulate shock waves and rarefaction waves 
for the Riemann problem of the system {\eqref{CL}}--{\eqref{PV}}. 
In Section {\ref{3}}, we find the admissible shock waves and rarefaction waves in settled regime by 
using the formula given in Section 2.  
\section{Preliminaries}
\label{2}
We let 
\begin{eqnarray*}
U=
\begin{pmatrix}
h  \\ n  
\end{pmatrix},  
\hspace{0.5cm}
F(U) = 
\begin{pmatrix}
\displaystyle\frac{1}{3}h^3   \\
\sqrt{\displaystyle\frac{2}{9C}}(nh)^{3/2}
\end{pmatrix},  
\end{eqnarray*}
so that the system {\eqref{CL}} and {\eqref{PV}} can be rewritten in the form 
\begin{eqnarray}
\partial_t U +\partial_x(F(U)) =0.  \label{system} 
\end{eqnarray}
It is well known that a solution to conservation laws {\eqref{system}} can become 
discontinuous even if the initial data is smooth. 
Therefore we treat the weak solution which is defined as follows : 
\begin{definition}[\cite{S}]
A bounded measurable function $U(x, t)$ is called a weak solution 
of the initial-value problem for {\eqref{system}} with bounded 
and measurable initial data $U(x, 0)$, provided that 
\begin{eqnarray}
\int^{\infty}_0 \int_{\R} (U\psi_t + F(U)\psi_x)dxdt 
+ \int_{\R} U(x, 0) \psi(x, 0) dx = 0 \label{weak}
\end{eqnarray}
holds for all $\psi \in C_0^1(\R\times\R_{+};\R^2)$. 
\end{definition} 
If the weak solution $U(x, t)$ has a discontinuity along a curve $x = x(t)$, 
the solution $U$ and the curve $x = x(t)$ must satisfy the 
Rankine-Hugoniot relations (jump conditions)
\begin{eqnarray}
s(U_L-U_R) = F(U_L)-F(U_R),  \label{RH1}
\end{eqnarray}
where $U_L=U(x(t)\,-\,0, t)$ is the limit of $U$ approaching $(x,t)$ from the
left and $U_R=U(x(t)\,+\,0, t)$ is the limit of $U$ approaching $(x,t)$ from 
the right, and $s = \frac{dx}{dt}$ is the propagation speed of $x(t)$. 
\par 
We consider the Riemann problem for the conservation laws {\eqref{system}} with the initial data 
called the Riemann data 
\begin{eqnarray}
U(x, 0)=
\begin{cases}
U_0 \hspace{0.5cm}x<0 \\
U_2 \hspace{0.5cm}x>0 
\end{cases}
. \label{initial}
\end{eqnarray} 
The Jacobian matrix of $F$ at $U$ is 
$$
DF(U) = 
\begin{pmatrix}
 h^2 & 0 \\
\displaystyle\sqrt{\frac{1}{2C}n^3 h } 
     &     
\displaystyle\sqrt{\frac{1}{2C}h^3 n}
\end{pmatrix}
$$ 
and district eigenvalues of $DF(U)$ are 
\begin{eqnarray}
\lambda_1 (U) = \sqrt{\frac{1}{2C}h^3 n}, 
\hspace{0.5cm}
\lambda_2 (U) =  h^2. \label{eigen}    
\end{eqnarray}
Here we assume that $h$ and $n$ are real valued function of $(x,t)\in \R\times \R_{+}$. According to {\cite{MPPB}}, set $C=2.307$ and $n=\phi h$,    
where the particle volume fraction $\phi$ satisfies $0\leq \phi< 1$.  
Under these conditions, the system {\eqref{system}} is strictly hyperbolic, 
i.e., district eigenvalues $\lambda_j(U)$ $(j=1,2)$ are real-valued and $\lambda_1(U)<\lambda_2(U)$ holds for any $U\in\Omega$, 
where 
$\Omega = \{(h,n)\in\R^2 : h> 0,\,\,0\leq n< h\}$. The right eigenvectors 
corresponding to the eigenvalues $\lambda_j(U)$ are 
\begin{eqnarray*}
r_1(U)= 
\begin{pmatrix}
0\\  
\vspace*{-0.3cm}\\
t_1 
\end{pmatrix}, 
\hspace{0.5cm}  
r_2(U)= 
\begin{pmatrix}
\displaystyle h^2- \sqrt{\frac{1}{2C}h^3 n} \\  
\vspace*{-0.3cm}\\
\displaystyle \sqrt{\frac{1}{2C}n^3 h } 
\end{pmatrix}, 
\end{eqnarray*}
where $t_1\neq 0$ is a constant. Note that  
$\nabla\lambda_1\cdot r_1=\frac{t_1}{2}\sqrt{\frac{1}{2Cn}h^3} \neq 0$ and 
$\nabla\lambda_2\cdot r_2=2h (h^2- \sqrt{\frac{1}{2C}h^3 n})\not= 0$ in $\Omega$, 
namely, the first and the second characteristic fields are genuinely nonlinear in $\Omega$. 
In this case, the weak solution will consist of 
three constant states $U_0,\,U_1,\,U_2$; the constant states $U_{j-1}$ and $U_j$ $(j=1,2)$ are connected by either shock waves or rarefaction waves 
(see {\cite{La1}, \cite{S}}). 
\par 
Fix the reference point $U_p=(h_p, n_p)$. 
We consider right states $U_R=U=(h, n)$ which can be connected 
to a left state $U_L=U_p$ followed by shock waves or rarefaction waves.  
If the weak solution has a jump discontinuity between the left state
$U_p$ and the right state $U$, then $U$ must satisfy 
the Rankine-Hugoniot relation {\eqref{RH1}}:  
\begin{eqnarray}
s(h-h_p) &=& \frac{1}{3}\left(h^3 - h_p^3\right), \label{RH}\\
s(n-n_p) &=& \sqrt{\frac{2}{9C}}
\left( (nh)^{3/2} - (n_ph_p)^{3/2}\right).   \nonumber  
\end{eqnarray}
Eliminating $s$ from these equations, we obtain 
\begin{eqnarray*}
(n-n_p)\left(h^2 + hh_p + h_p^2 \right) = \sqrt{\frac{2}{C}}
\left( (nh)^{3/2} - (n_ph_p)^{3/2}\right) \label{curve}
\end{eqnarray*}
whose graph is called the {\em Hugoniot locus}. 
In order to pick up physically relevant solutions, we further require the following 
$k$-entropy inequalities $(k=1,2)$   
\begin{eqnarray*}
&&s<\lambda_1(U_p), \hspace{0.5cm}\lambda_1(U)< s< \lambda_2(U), 
\hspace{0.5cm} \text{($1$-entropy inequality)},
\label{1entropy} \\
&& \nonumber \\
&&\lambda_1(U_p)< s< \lambda_2(U_p),\hspace{0.5cm}\lambda_2(U)<s, 
\hspace{0.5cm} \text{($2$-entropy inequality)},  
\label{2entropy}
\end{eqnarray*}
which in this case reads 
\begin{eqnarray}
&&\sqrt{\frac{1}{2C}h^3 n} <s < 
\min \bigg\{\sqrt{\frac{1}{2C}h_p^3 n_p},\,h^2 \bigg\}, \hspace{0.5cm} 
\text{($1$-entropy inequality)},  \label{1entro}   \\
&&                                        \nonumber \\
&&\max \bigg\{\sqrt{\frac{1}{2C}h_p^3 n_p},\,h^2 \bigg\} <s < h^2_p, 
\hspace{0.5cm} 
\text{($2$-entropy inequality)},  \label{2entro}
\end{eqnarray}
where $s$ is the speed of discontinuity
\begin{eqnarray*}
s= \left(\frac{2}{81C}\right)^{1/4}\sqrt{\frac{(h^2 + hh_p + h_p^2)
\left( (nh)^{3/2} - (n_ph_p)^{3/2}\right)}{n-n_p}}.  
\end{eqnarray*}
If $U$ satisfies {\eqref{RH}} and {\eqref{1entro}}, then 
$U$ can be connected to $U_p$ from the right followed by a {\em 1-shock wave}.
Since the system {\eqref{system}} is strictly hyperbolic, it is clear that 
$\sqrt{\frac{1}{2C}h^3 n}<h^2$. Thus the $1$-shock curve is given by 
\begin{align}
\notag
S_1(U_p) = \{(h, n) : &(n-n_p)\left(h^2 + hh_p + h_p^2 \right) \\
\label{1shock} 
&= \sqrt{\frac{2}{C}}\left( (nh)^{3/2} - (n_ph_p)^{3/2}\right),
h^3n < h^3_p n_p \}. 
\end{align}
Similarly, $U$ can be connected to 
$U_p$ from the right followed by a {\em 2-shock wave}, 
provided $U$ satisfies {\eqref{RH1}} and {\eqref{2entro}}.    
This curve is called the $2$-shock curve, which is given by 
\begin{align}
\notag
S_2(U_p) = \{(h, n) : &(n-n_p)\left(h^2 + hh_p + h_p^2 \right) \\
&= \sqrt{\frac{2}{C}}\left( (nh)^{3/2} - (n_ph_p)^{3/2}\right),
h<h_p\}. 
\label{2shock}   
\end{align}
\par
We consider candidates of right states $U_R=U=(h, n)$ which can 
be connected to a given left state $U_L=U_p=(h_p, n_p)$ followed by  
a {\em rarefaction wave}. Here we note that the condition for (physically relevant)
rarefaction waves is that the corresponding eigenvalue (speed) 
$\lambda$ increases from the left to the right side of the wave 
(see {\cite{S}}), that is  
\begin{eqnarray}
\lambda (U_p) < \lambda (U).  \label{speedrare}
\end{eqnarray}
The Riemann problem {\eqref{system}}, {\eqref{initial}} are invariant 
under the scaling $(t, x)\mapsto(\eta t, \eta x)$ for all $\eta >0$. 
Therefore we seek self-similar solutions of the form 
$U(x,t)\equiv U(\frac{x}{t})$. 
If we let $\xi=\frac{x}{t}$, then we see that $U(\xi)$ satisfies the ordinary differential equation 
\begin{eqnarray*}
(DF(U)-\xi) d_{\xi}U=0,   
\end{eqnarray*}
where $d_{\xi}=\frac{d}{d\xi}$.  
If $d_{\xi}U\neq 0$, then $\xi$ is the eigenvalue for $DF(U)$ and $d_{\xi}U$ is the corresponding eigenvector. 
Since $DF(U)$ has two real and distinct eigenvalues $\lambda_1<\lambda_2$, there exist two families of rarefaction waves, 
{\em $1$-rarefaction waves} and {\em $2$-rarefaction waves}.
For $1$-rarefaction waves, the eigenvector $d_{\xi}U=(d_{\xi}h, d_{\xi}n)^{\top}$ satisfies 
\begin{eqnarray*}
(-\lambda_1(U)\I+DF(U)) d_{\xi}U
= 
\begin{pmatrix}
-\displaystyle\sqrt{\frac{1}{2C}h^3 n}+h^2 & 0 \\
\displaystyle\sqrt{\frac{1}{2C}n^3 h } 
     &     
0
\end{pmatrix}
\begin{pmatrix}
 d_{\xi}h \\ 
 d_{\xi}n 
\end{pmatrix}
=
\begin{pmatrix}
 0 \\ 
 0
\end{pmatrix}
,
\end{eqnarray*}
which gives $d_{\xi}h=0$. Since $d_{\xi}n\neq 0$, we have 
\begin{eqnarray*}
\displaystyle
\frac{dh}{dn}= 0. 
\end{eqnarray*}
We integrate this to obtain the curve passing all possible $U$ connected to 
$U_p$ followed by a $1$-rarefaction wave. This curve is called the $1$-rarefaction 
curve, which is in our case given by 
\begin{eqnarray}
R_1(U_p) = \{(h, n) : h=h_p,\,\,n>n_p \},      
\label{1rare} 
\end{eqnarray}
where $n>n_p$ comes from $\lambda_1(U_p)<\lambda_1(U)$.
\par 
For $2$-rarefaction waves, the eigenvector $d_{\xi}U$ satisfies 
\begin{eqnarray*}
(-\lambda_2(U)\I+DF(U)) d_{\xi}U
= 
\begin{pmatrix}
0 & 0 \\
\displaystyle\sqrt{\frac{1}{2C}n^3 h } 
     &     
-h^2+\displaystyle\sqrt{\frac{1}{2C}h^3 n}
\end{pmatrix}
\begin{pmatrix}
 d_{\xi}h \\ 
 d_{\xi}n 
\end{pmatrix}
=
\begin{pmatrix}
 0 \\ 
 0
\end{pmatrix}
,
\end{eqnarray*}
which gives 
\begin{eqnarray*}
\displaystyle
\frac{dh}{dn}= \frac{h^2- \sqrt{\frac{1}{2C}h^3 n}}
{\sqrt{\frac{1}{2C}n^3 h}}= \left(\sqrt{2C}\sqrt{\frac{h}{n}}-1\right)\frac{h}{n}. 
\end{eqnarray*}
We can solve this ordinary differential equation, the solution is given by  
\begin{eqnarray*}
h=\frac{n}{(\sqrt{\frac{C}{2}} -e^{A}n)^2}, 
\end{eqnarray*}
where $e^{A}$ is the constant of integration. 
When the solution takes $U_p=(h_p, n_p)$, the constant $e^A$ is determined as $\frac{1}{n_p}(\sqrt{\frac{C}{2}}-\sqrt{\frac{n_p}{h_p}})$  
then the special solution is obtained as 
\begin{eqnarray*}
h=\frac{n\,n_p^2}{\left(n\sqrt{\frac{n_p}{h_p}}-(n-n_p)\sqrt{\frac{C}{2}} 
\right)^2}. 
\end{eqnarray*}
The graph of this function is called the $2$-rarefaction curve consisting of $U$ which can be connected from the left state $U_p$ by a 
$2$-rarefaction wave. 
We denote by 
\begin{eqnarray}
R_2(U_p) = \{(h, n) : h\left(n\sqrt{\frac{n_p}{h_p}}-(n-n_p)
\sqrt{\frac{C}{2}}\right)^2=n\,n_p^2,\,\,h >h_p \}, 
\label{2rare} 
\end{eqnarray}
where the condition $h_p < h$ comes from $\lambda_2(U_p)<\lambda_2(U)$.

\section{Admissble weak solutions for the settled regime}
\label{3}
In this section we construct weak solutions of Riemann problem 
{\eqref{system}}, {\eqref{initial}} by substituting the values corresponding to 
the settle regime into the curves given in the previous section. 
We tackle the Riemann problem for situations wherein $h<h_p$ and $h>h_p$ 
representing a step-down and step-up function, respectively. 
\par 
We begin with finding admissible wave curves connecting from the fixed left state 
$U_0$ to the right states $U=(h, n)$ when $h < h_0$.  
We set $U_0=(h_0, n_0)=(1, 0.1)$ and $C=2.307$, which are used in {\cite{MPPB}}. 
Then the Hugoniot locus becomes the set  
\begin{eqnarray}
S(U_0) : \left\{(n-0.1)\left(h^2 + h + 1 \right) = \sqrt{\frac{2}{2.307}}
\left( (nh)^{3/2} - (0.1)^{3/2}\right)\right\}, \label{curve1}
\end{eqnarray}
while the $1$-entropy inequality and the $2$-entropy inequality are as follows, 
respectively :  
\begin{eqnarray}
\sqrt{\frac{1}{4.614}h^3 n} <
&s&
< \min \bigg\{\sqrt{\frac{1}{46.14}},\,h^2 \bigg\},  \label{1entro1} \\
&&\nonumber \\
\max \bigg\{\sqrt{\frac{1}{46.14}},\,h^2 \bigg\} <
&s& 
< 1, \label{2entro1}
\end{eqnarray}
where 
\begin{eqnarray}
s=\displaystyle\left(\frac{2}{186.867}\right)^{1/4}\sqrt{\frac{(h^2 + h + 1)\left( (nh)^{3/2} - (0.1)^{3/2}\right)}{n-0.1}}.  \label{speed}
\end{eqnarray}
We note that inequalities {\eqref{1entro1}}, {\eqref{2entro1}} are equivalent to the following inequalities : 
\begin{eqnarray}
s-\sqrt{\frac{1}{4.614}h^3 n}>0
&\hspace{0.5cm}\hbox{and}\hspace{0.5cm}&
s-\min \bigg\{\sqrt{\frac{1}{46.14}},\,h^2\bigg\} 
< 0,  \label{bound}\\
s-\max \bigg\{\sqrt{\frac{1}{46.14}},\,h^2 \bigg\} >0
&\hspace{0.5cm}\hbox{and}\hspace{0.5cm}&
s-1 < 0. \label{bound2}
\end{eqnarray}
\begin{figure}[h]
    \begin{center}
    \subfloat[$\lambda_1(U_0)\geq \lambda_2(U)$ ]{
            \includegraphics[scale=0.6]{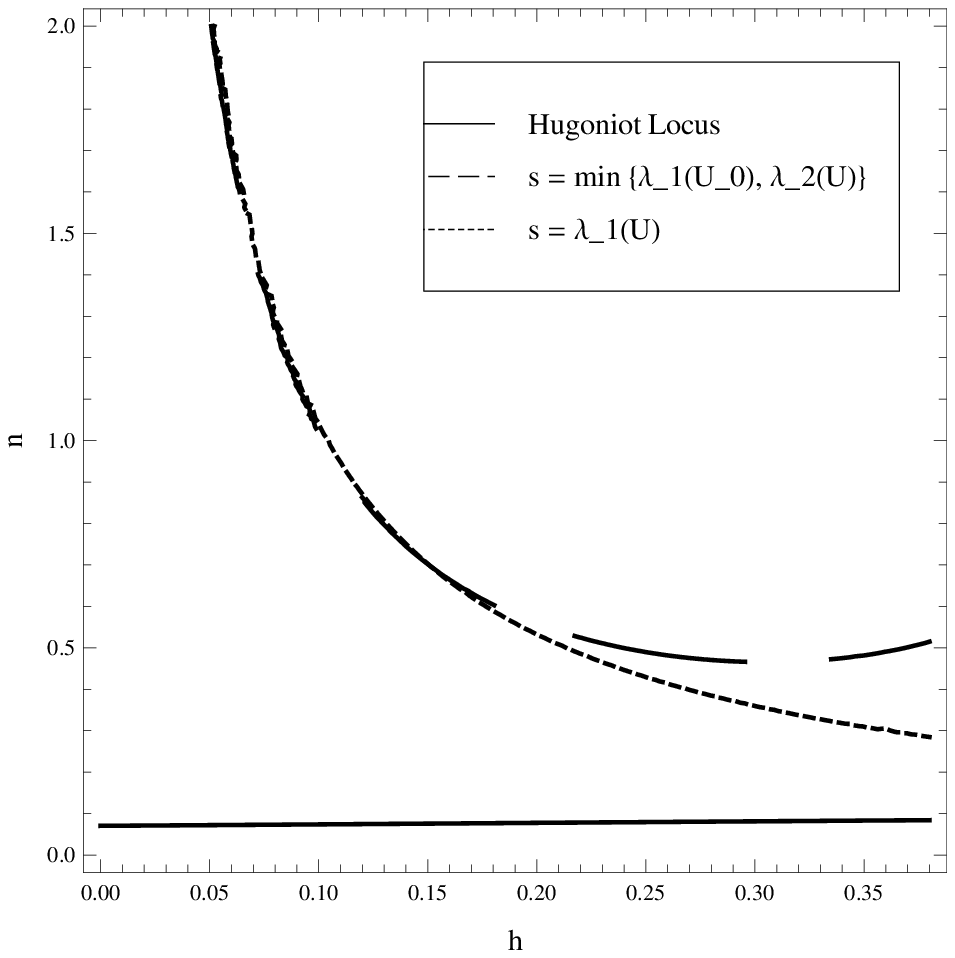}
\label{1s1}}
\subfloat[$\lambda_1(U_0)<\lambda_2(U)$]
{\includegraphics[scale=0.6]{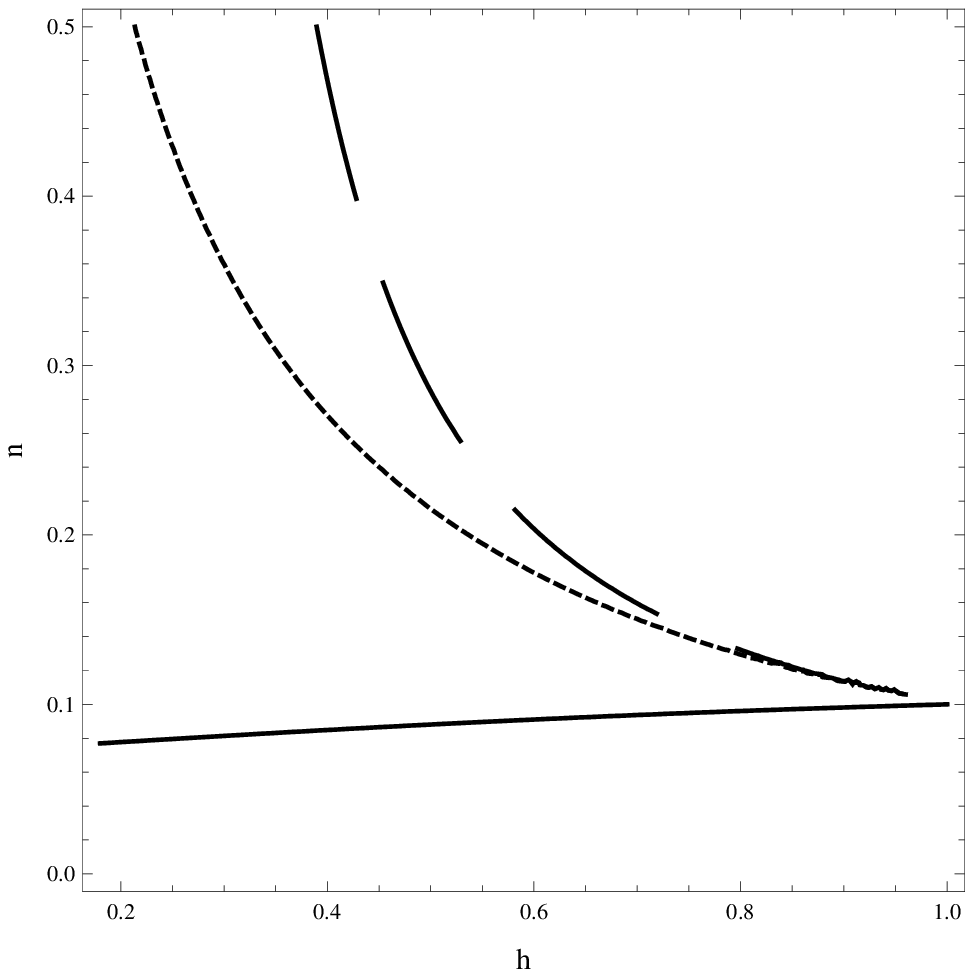}
\label{1s2}}
\caption{Hugoniot locus and the $1$-entropy inequality. We plot the Hugoniot locus {\eqref{curve1}} and implicit functions $s=\lambda_1(U)$ 
and $s=\min\{\lambda_1(U_0), \lambda_2(U)\}$, where $\lambda_1(U)=\sqrt{\frac{1}{4.614}h^3 n}$, 
$\lambda_1(U_0)=\sqrt{\frac{1}{46.14}}$, $\lambda_2(U)=h^2$. The solid, dashed and dotted curves represent 
the Hugoniot locus {\eqref{curve1}}, $s=\min\{\lambda_1(U_0), \lambda_2(U)\}$ and $s=\lambda_1(U)$ respectively.}
{\label{1shockw}}
    \end{center}    
\end{figure}
\par 
We shall examine whether there exists $(h, n)$ satisfying 
{\eqref{curve1}} and {\eqref{bound}} with phase portraits. 
In Figure {\ref{1shockw}}
we plot the Hugoniot locus {\eqref{curve1}} and the implicit functions 
$s=\sqrt{\frac{1}{4.614}h^3 n}$ and 
$s=\min \{\sqrt{\frac{1}{46.14}},\,h^2 \}$, which is $s=h^2$ 
for the case $\sqrt{\frac{1}{46.14}}\geq h^2$ 
(Figure {\ref{1shockw}\subref{1s1}}) and 
$s=\sqrt{\frac{1}{46.14}}$ 
for the case $\sqrt{\frac{1}{46.14}}< h^2$ 
(Figure {\ref{1shockw}\subref{1s2}}). 
Two dashed lines in Figure {\ref{1shockw}} show the upper bound and lower bound 
for the inequality {\eqref{bound}}, which means that every point $(h, n)$ within 
the open region between the upper graph $s=\min \{\sqrt{\frac{1}{46.14}},\,h^2 \}$ 
and the lower graph $s=\sqrt{\frac{1}{4.614}h^3 n}$ satisfies {\eqref{bound}}. 
As can be seen from the figure, $(h, n)$ satisfying the Rankine-Hugoniot relation {\eqref{curve1}} 
does not belong to the region that the 
$1$-entropy inequality {\eqref{bound}} holds. 
Thus, the weak solution does not admit 1-shock waves.  
\begin{figure}[h]
\begin{center}
    \subfloat[$\lambda_1(U_0)\geq \lambda_2(U)$ ]{
            \includegraphics[scale=0.6]{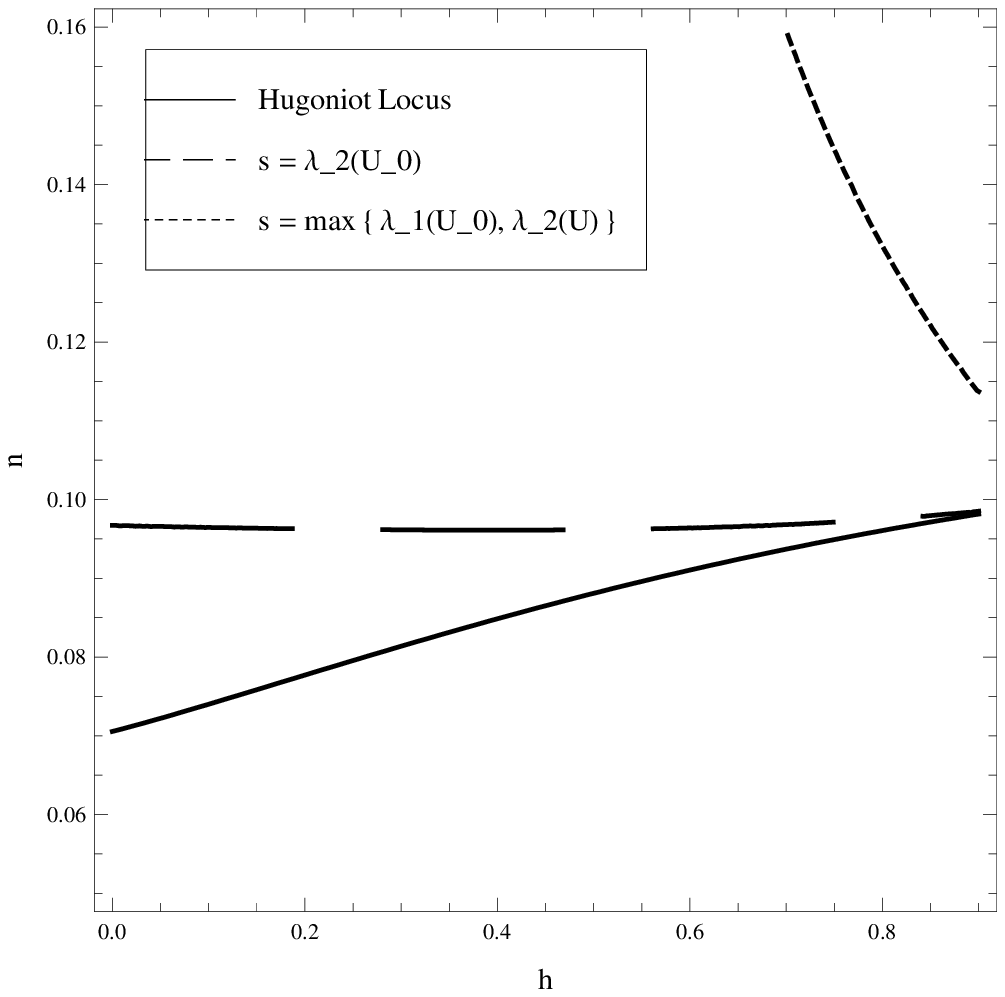}
\label{2s1}}
\subfloat[$\lambda_1(U_0)<\lambda_2(U)$]
{\includegraphics[scale=0.6]{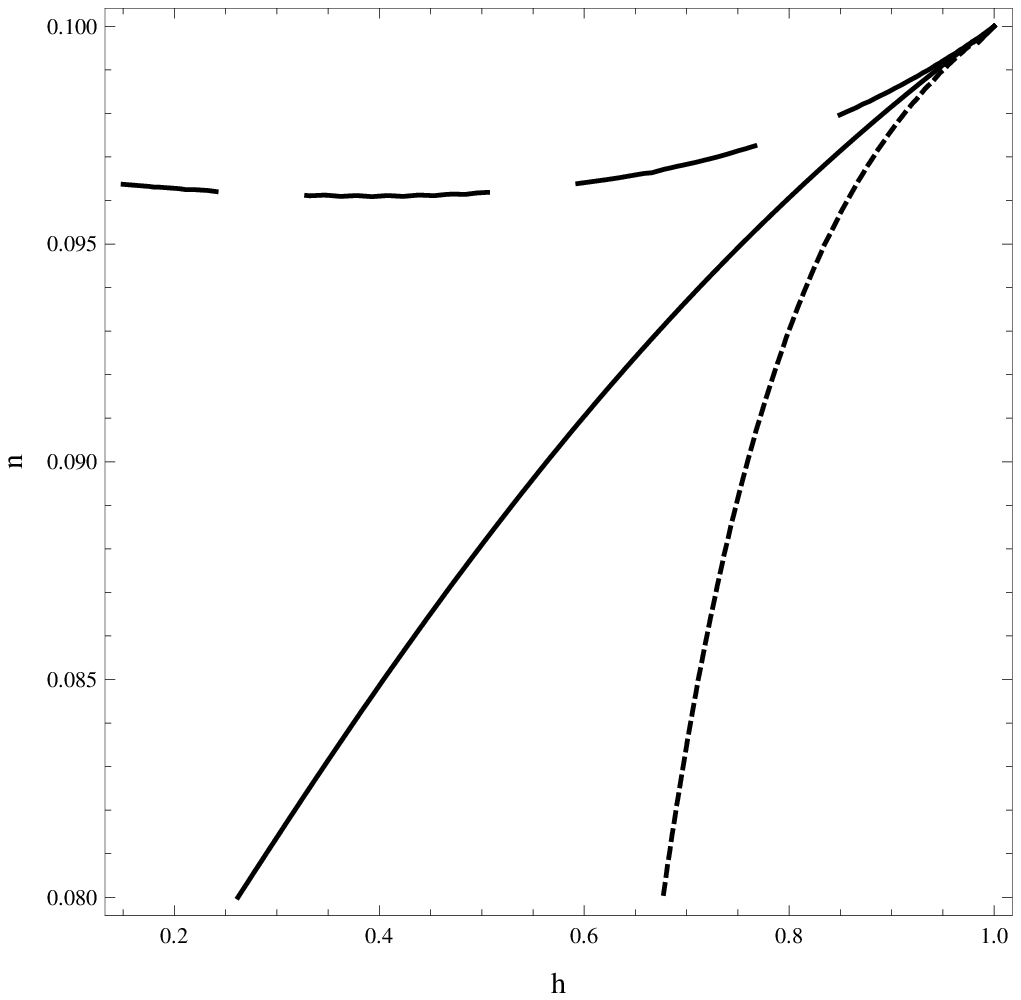}
\label{2s2}}
\caption{Hugoniot locus and the $2$-entropy inequality. In this figure we plot the Hugoniot locus {\eqref{curve1}} and implicit functions $s=\lambda_2(U_0)$ 
and $s=\max\{\lambda_1(U_0), \lambda_2(U)\}$, where $\lambda_1(U_0)=\sqrt{\frac{1}{46.14}}$, 
$\lambda_2(U_0)=1$, $\lambda_2(U)=h^2$. The solid, dashed and dotted lines represent the Hugoniot locus {\eqref{curve1}}, 
$s=\lambda_2(U_0)$ and $s=\max\{\lambda_1(U_0), \lambda_2(U)\}$ respectively.}
{\label{2shockw}} 
\end{center}
\end{figure}
\par 
Similarly, we examine whether there exists a (right) state $(h, n)$ satisfying 
{\eqref{curve1}} and {\eqref{bound2}}. 
In Figure {\ref{2shockw}} we plot the Hugoniot locus {\eqref{curve1}} 
and the implicit functions $s=1$ 
and $s=\max \{\sqrt{\frac{1}{46.14}},\,h^2 \}$. 
When $\sqrt{\frac{1}{46.14}}\geq h^2$, every point $(h, n)$ 
satisfying the Rankine-Hugoniot relation {\eqref{curve1}} does not belong to the region 
between the upper graph $s=1$ and the lower graph  
$s=h^2$ (Figure {\ref{2shockw}\subref{2s1}}). On the other hand, when 
$\sqrt{\frac{1}{46.14}}< h^2$, the Hugoniot locus $S(U_0)$ belongs 
to the region between the upper graph $s=1$ and the lower graph 
$s=\sqrt{\frac{1}{46.14}}$ (Figure {\ref{2shockw}\subref{2s2}}), which means that 
there exists $(h, n)$ satisfying both {\eqref{curve1}} and 
{\eqref{bound2}}. Thus, when $\sqrt{\frac{1}{46.14}}< h^2$, the $2$-shock 
wave exists and the $2$-shock curve is given by {\eqref{curve1}} for $h<1$.  
\begin{figure}[h]
    \begin{center}
        {\includegraphics[scale=0.65]{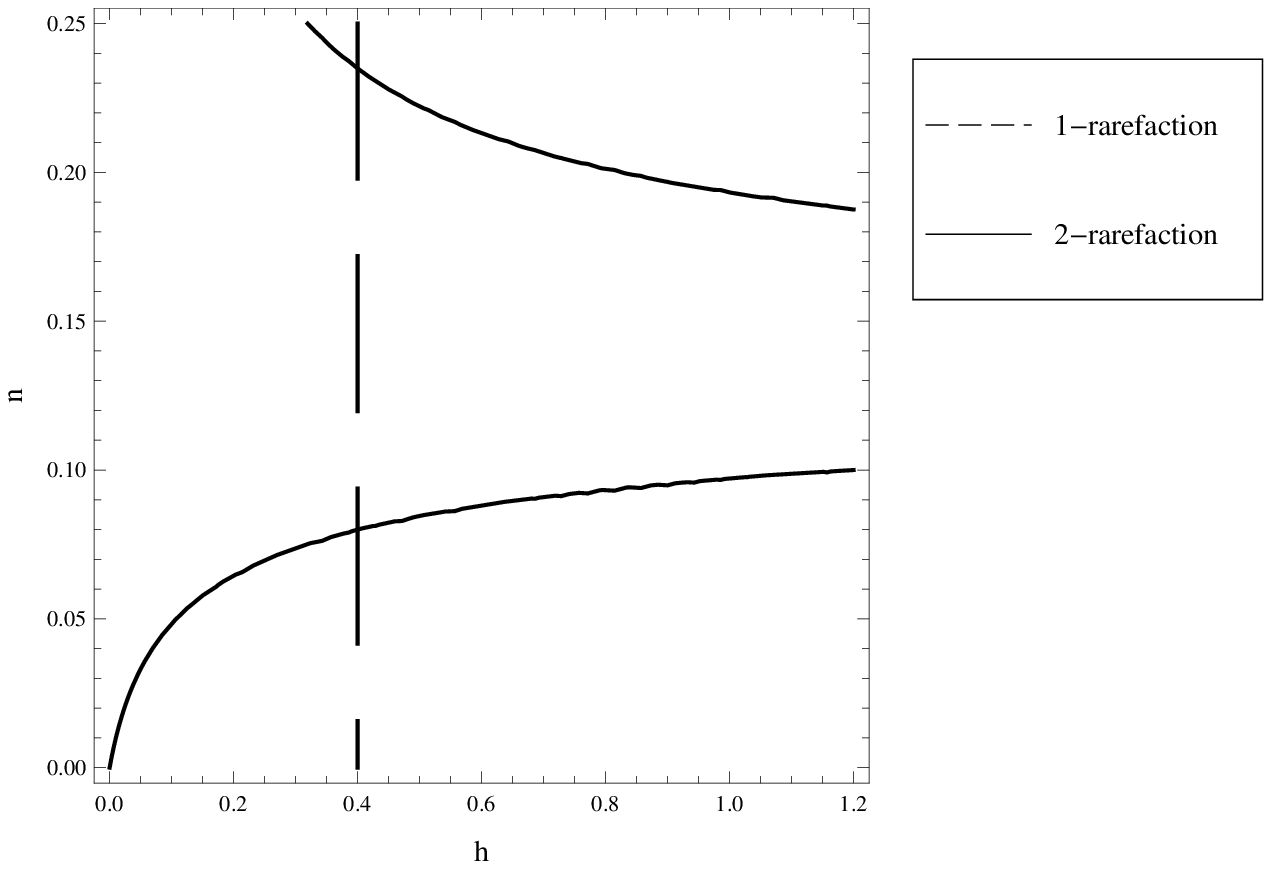}}
\caption{In this figure we plot a graph of two rarefaction wave curves {\eqref{1rare2}} and {\eqref{2rare2}}. 
The dashed and solid lines represent the $1$-rarefaction wave curve {\eqref{1rare2}} and 
the $2$-rarefaction wave curve {\eqref{2rare2}} respectively.}
{\label{weak1}}
    \end{center}    
\end{figure}
\par 
As a example, we take $U_2=(0.2, n_{2,s})$
{\footnote{Using Newton's method, a sample of the approximate solution for equation {\eqref{eqs}} is obtained as $n_{2,s}=0.0777100325$. }}, 
where $n_{2,s}$ is the solution of 
\begin{eqnarray}
1.24\,(n_{2,s}-0.1) = \sqrt{\frac{2}{2.307}}
\left( (0.2\,n_{2,s})^{3/2} - (0.1)^{3/2}\right), \label{eqs}
\end{eqnarray}
which is exactly the equation {\eqref{curve1}} with $U=U_2$.
Then the left state $U_0=(1, 0.1)$ and the right state $U_2$ is connected by a single $2$-shock wave.  
In the range $h<1$,{\eqref{1shock}} and {\eqref{2shock}} 
make no sense as $1$-shock wave and $2$-shock wave by the entropy inequalities, respectively.  
\bigskip
\par 
Similarly, we find admissible wave curves in the case $h > h_0$.  
Fix $U_0=(h_0, n_0)=(0.4, 0.08)$ and $C=2.307$, and  
we plot the $1$-rarefaction curve and $2$-rarefaction curve, which are given 
as follows, respectively :  
\begin{eqnarray}
&h= 0.4,&\hspace{0.5cm}n>0.08,    \label{1rare2} \\
&h= \frac{n\,(0.08)^2}{\left(n\sqrt{0.2}-(n-0.08)
\sqrt{\frac{2.307}{2}}\right)^2}, &\hspace{0.5cm}
h>0.4,  \label{2rare2}
\end{eqnarray}
which means that {\eqref{1rare2}} makes no sense as $1$-rarefaction
{\footnote{When $h\not = h_p$, which is typical as phenomena of fluid motion {\cite{MB}}, $1$-rarefaction waves do not exist. 
On the other hand, if we admit 
$h=h_p$, a $1$-rarefaction wave connecting $(h_p,n_p)$ and $(h_p,n)$ with $n_p < n < h_p$ is also admitted. }, 
but {\eqref{2rare2}} makes sense as $2$-rarefaction by {\eqref{speedrare}}. 
\par 
As an example, we take $U_2=(1.0, n_{2,r})$ 
{\footnote{Using Newton's method, a sample of the approximate solution for equation {\eqref{eqr}} is obtained as $n_{2,r}=0.0972723141$. }},  
where $n_{2,r}$ is the solution of 
\begin{eqnarray}
0.08\sqrt{n_{2,r}}+(n_{2,r}-0.08)\sqrt{\frac{2.307}{2}} = n_{2,r}\sqrt{0.2}. \label{eqr}
\end{eqnarray}
Then the left state $U_0=(0.4, 0.08)$ and the right state $U_2$ is connected by a single $2$-rarefaction wave.  
\par 
Our argument is summarized in Table {\ref{list}}. 
Following the terminology \lq\lq allowed sequence" of waves in \cite{SMP}, wave sequences consisting of shocks and rarefactions associated with 
{\em identical} characteristic fields are excluded.
\begin{table}
\begin{tabular}{|c|c|c|}
\hline
$w_1$&$w_2$&appear \\  
$   $&$   $&ance \\ 
\hline
\multicolumn{2}{|c|}{$1$-rarefaction}&$\triangle$\\
\hline
$1$-rarefaction &$2$-rarefaction&$\triangle$
\\
\hline
$1$-rarefaction&$1$-shock wave &$\times$\\
\hline
$1$-rarefaction & $2$-shock wave &$\times$\\
\hline\hline 
$2$-rarefaction &$1$-rarefaction &$\times$\\
\hline
\multicolumn{2}{|c|}{$2$-rarefaction}&$\bigcirc$\\
\hline
$2$-rarefaction &$1$-shock wave &$\times$\\
\hline
$2$-rarefaction &$2$-shock wave &$\times$\\
\hline
\end{tabular}
\begin{tabular}{|c|c|c|}
\hline
$w_1$&$w_2$&appear\\ 
$   $&$   $&ance \\ 
\hline
$1$-shock wave&$1$-rarefaction&$\times$\\
\hline
$1$-shock wave&$2$-rarefaction&$\times$\\
\hline
\multicolumn{2}{|c|}{$1$-shock wave}&$\times$\\
\hline
$1$-shock wave&$2$-shock wave&$\times$\\
\hline\hline 
$2$-shock wave &$1$-rarefaction &$\times$\\
\hline
$2$-shock wave &$2$-rarefaction &$\times$\\
\hline
$2$-shock wave &$1$-shock wave &$\times$\\
\hline
\multicolumn{2}{|c|}{$2$-shock wave}&$\bigcirc$\\
\hline
\end{tabular} 
\caption{Combination of solutions to appearance. $w_i$ $(i=1,2)$ denote the simple wave in the $i$-characteristic field.}
\label{list}
\end{table}

\section{Conclusions}
\label{5}
In this paper we have dealt with a Riemann problem for the system of conservation laws {\eqref{CL}}--{\eqref{PV}} 
which is derived from the dilution approximation of a suspension flow on an incline as a mathematical model in the settled regime. 
Murisic et al. {\cite{MPPB}} dealt only with a exact solution for the system {\eqref{CL}}--{\eqref{PV}}, 
when the initial volume fraction is fixed as $\phi(0, x)\equiv f_0$ for some given $f_0\ll1$.     
On the other hand, we aim at covering the solution of this system when the initial volume fraction $\phi(0,x)$ is a variable satisfying $0<\phi<1$.     
In Sections {\ref{2}} and {\ref{3}}, we show that the weak solution of this Riemann problem is connected by a single $2$-rarefaction wave 
from the left state $U_0=(h_0, n_0)$ to the right state $U_2=(h_2, n_2)$ when $h_0<h_2$, and connected by a single $2$-shock wave when $h_0>h_2$. 
To illustrate one example of these wave curves, we impose the initial conditions as follows,    
\begin{eqnarray*}
U^r(x, 0)=
\begin{cases}
U_0=(0.4, 0.08)    &x<0 \\
U_2=(1.0, n_{2,r}) &x>0 
\end{cases}
, 
\hspace{0.8cm}
U^s(x, 0)=
\begin{cases}
U_0=(1.0, 0.1)     &x<0 \\
U_2=(0.2, n_{2,s}) &x>0 
\end{cases}
,
\end{eqnarray*}
where $n_{2, s}$ and $n_{2, r}$ is the solution of {\eqref{eqs}} and {\eqref{eqr}} respectively. 
We take the values of $U^r(x, 0)$ and $U^s(x, 0)$ to satisfy the ranges $0\leq h\leq 1$ and $0\leq n\leq 0.1$ of the exact solution 
handled in {\cite{MPPB}}.
With the Riemann data $U^r(x,0)$, the weak solution consists of a single $2$-rarefaction wave whose curve is shown in Figure {\ref{weak1}}. 
With the Riemann data $U^s(x,0)$, the weak solution consists of a single $2$-shock wave whose curve is shown in Figure {\ref{2shockw}\subref{2s2}}.  
The construction method given in Sections {\ref{2}} and {\ref{3}} may also be useful for other suspension models even if the initial volume fraction 
$\phi$ depends on $(t,x)$.
In the future work, we will investigate whether rarefaction wave and shock wave obtained from {\eqref{CL}} and {\eqref{PV}} 
correspond to experiment results in {\cite{MPPB}}.

\section*{Acknowledgements}
\label{ac}
This work is supported by $2017$ IMI Joint Use Research Program CATEGORY 
\lq\lq Short-term Visiting Researcher\rq\rq  in Institute of Mathematics for Industry, Kyushu University. 
KM was partially supported by Program for Promoting the reform of national universities (Kyushu University), 
Ministry of Education, Culture, Sports, Science and Technology (MEXT), Japan, World Premier International Research Center Initiative (WPI), MEXT, Japan 
and JSPS Grant-in-Aid for Young Scientists (B) (No. JP17K14235). 



\begin{thebibliography}{10}
\bibitem{H1}H. Huppert. Flow and instability of a viscous current down a slope. Nature, 
{\bf {300}} 427--429, (1982).
\bibitem{La1}P. D. Lax, Hyperbolic system of conservation laws II, 
Comm. Pure Appl. Math. {\bf{10}} 537--566, (1957). 
\bibitem{MB}A. Mavromoustaki, A. L. Bertozzi, Hyperbolic systems of conservation laws in gravity-driven, particle-laden thin-film flows, 
Journal of Engineering Mathematics {\bf{88}} 29--48, (2014).
\bibitem{MPPB}N. Murisic, B. Pausader, D. Peschka, A. L. Bertozzi, Dynamics of particle settling and resuspension in viscous liquids, 
J. Fluid Mech. {\bf{717}} 203--231, (2013). 
\bibitem{SMP}S. Schecter, D. Marchesin, B. J. Plohr, 
Structurally stable Riemann solutions, 
J. Differential Equations {\bf{126}} no. 2, 303--354, (1996). 
\bibitem{S}J. Smoller. Shock waves and reaction-diffusion equations, 
{\bf 258} of Grundlehren der Mathematischen Wissenschaften (Fundamental Principles of Mathematical Sciences). Springer-Verlag, New York, 
second edition, (1994). 
\end{thebibliography}
\end{document}